\documentclass[conference]{IEEEtran}
\usepackage{amssymb}
\usepackage{amsfonts}
\usepackage{amsmath}

\newtheorem{remark}{Remark}

\DeclareMathOperator{\Ker}{Ker}

\begin{document}
\title{Pseudo prolate spheroidal functions}

\author{\IEEEauthorblockN{Lu\'{\i}s Daniel Abreu}
\IEEEauthorblockA{Austrian Academy of Sciences\\
Acoustics Research Institute \\
Wohllebengasse 12-14\\
Vienna A-1040, Austria\\
Email: daniel@mat.uc.pt}
\and
\IEEEauthorblockN{Jo\~{a}o M. Pereira}
\IEEEauthorblockA{Program in Applied and\\
Computational Mathematics\\
Princeton University \\
Princeton, NJ 08544, USA\\
Email: jpereira@princeton.edu}}
\thanks{L. D.~Abreu was supported by Austrian Science Foundation (FWF)
START-project FLAME ("Frames and Linear Operators for Acoustical Modeling
and Parameter Estimation", Y 551-N13).\\
J. M. Pereira was supported by...}

\maketitle

\begin{abstract}
Let $D_{T}$\ and $B_{\Omega }$\ denote the operators which cut the time
content outside $T$ and the\ frequency content outside $\Omega $,
respectively. The prolate spheroidal functions are the eigenfunctions of the
operator $P_{T,\Omega }=D_{T}B_{\Omega }D_{T}$. With the aim of formulating
in precise mathematical terms the notion of Nyquist rate, Landau and Pollack
have shown that, asymptotically, the number of such functions with
eigenvalue close to one is $\approx \frac{\left\vert T\right\vert \left\vert
\Omega \right\vert }{2\pi }$.\ We have recently revisited this problem with
a new approach: instead of counting the number of \emph{eigenfunctions with
eigenvalue close to one}, we count the maximum number of \emph{orthogonal }$%
\epsilon $\emph{-pseudoeigenfunctions with }$\epsilon $\emph{%
-pseudoeigenvalue one}. Precisely, we count how many orthogonal functions
have a maximum of energy $\epsilon $ outside the domain $T\times \Omega $,
in the sense that $\left\Vert P_{T,\Omega }f-f\right\Vert ^{2}\leq \epsilon $%
. We have recently discovered that the sharp asymptotic number is $\approx
(1-\epsilon )^{-1}\frac{\left\vert T\right\vert \left\vert \Omega
\right\vert }{2\pi }$. The proof involves an explicit construction of the
pseudoeigenfunctions of $P_{T,\Omega }$. When $T$ and $\Omega $ are
intervals we call them \emph{pseudo prolate spheroidal functions}. In this
paper we explain how they are constructed.
\end{abstract}

\maketitle

\section{Introduction}

\subsection{Slepian's bandwith paradox}

In his 1974 Shannon lecture, whose written version appeared in \cite{Slepian}%
, David Slepian stated the following paradoxal dilemma:

\emph{"It is easy to argue that real signals must be
bandlimited. It is also easy to argue that they cannot be
so."}

Such a dilemma (we will call it the \emph{bandwith paradox}) reflects a
mainstay of quantitative physical sciences: the gap between observations and
models of the real world. On the one side, it is reasonable to accept that,
for any measuring instrument, there is a finite cutoff above which the
instrument would not be able to measure the frequencies of a signal. Hence,
it can be argued that \emph{all signals are bandlimited}. On the other side,
bandlimited signals are represented by analytic functions. This does not
allow the function to vanish in any real interval, leading to the
unrealistic model where signals cannot start or stop, but must go on
forever. Hence, it can be argued that \emph{no signal is bandlimited}.

The heuristics of the previous paragraph are already enough to change our
mindset: instead of \emph{supports} one should think of\emph{\ essential
supports}. Then the question arises of what is the \emph{dimension} of the
set of such functions. Since in reality we are not dealing with finite
dimensonal sets, we need to resort to an \emph{approximated notion of
dimension}. For instance, Landau and Pollack \cite{LP} considered, as a
notion of dimension, the minimal number $N(\epsilon )$ of functions required
to approximate any essentially time-band limited function in the $L_{2}$
norm up to an error $\epsilon $. Based on such considerations, two solutions
of the bandwith paradox have been offered, one by Landau and Pollack, the
other by Slepian. We will give a brief acount of the two approaches and
suggest a new one, based on a line of research initiated in \cite{Measures}.

We note in passing that, besides the solution of the bandwith paradox, some
of the above heuristics played a fundamental role in the papers \cite%
{DonohoLogan} and \cite{DonohoStark}, which spearheaded the modern theory of
Compressed Sensing, where an understanding of the deep mathematical reasons
behind the sparsity-promoting properties of $l_{1}$ minimization has been
achieved\ \cite{Cand2014}.

\subsection{ Landau-Pollack solution: prolate spheroidal functions}

Let $D_{T}$\ and $B_{\Omega }$\ denote the operators which cut the time
content outside $T$ and the\ frequency content outside $\Omega $,
respectively.\ In a nowadays classical paper \cite{LP}, whose purpose was to 
\emph{examine the true in the engineering intuition that there are
approximately }$\left\vert T\right\vert \left\vert \Omega \right\vert /2\pi $%
\emph{\ independent signals of bandwidth }$\Omega $\emph{\ concentrated on
an interval of length }$T$, Landau and Pollak have considered the eigenvalue
problem associated with the positive self-adjoint operator $P_{T,\Omega
}=D_{T}B_{\Omega }D_{T}$. When $T$ and $\Omega $ are real intervals, $%
P_{T,\Omega }$ can be written explicitly as%
$$(P_{T,\Omega }f)(x)=\left\{ 
\begin{array}{c}
\int_{T}\frac{\sin \Omega (x-t)}{\pi (x-t)}f(t)dt\mbox{ \ \ \ }if\mbox{ \ }%
x\in T \\ 
0\mbox{ \ \ \ \ \ \ \ \ \ \ \ \ \ \ \ \ \ \ \ \ \ \ \ }if\mbox{ \ }x\notin T%
\end{array}%
\right. .
$$%
The eigenfunctions of $P_{T,\Omega }$ are the \emph{prolate spheroidal
functions} $\{\phi _{j}\}_{j=0}^{\infty }$.\ They provide the best known
dictionary for approximating essentially time and band limited functions in
the $L_{2}$ norm \cite{LP} and their properties are still object of current
investigation \cite{prolates}. The approach to the bandwith paradox based on
prolate spheroidal functions relies on the peculiar behaviour of the spectra
of $P_{T,\Omega }$: the largest eigenvalues of $P_{T,\Omega }$ are very
close to $1$, before plunging very fast to almost $0$. But the eigenvalues
of $P_{T,\Omega }$ are the singular values of the operator $B_{\Omega }D_{T}$%
, whose singular functions satisfy 
\begin{equation}
\int_{_{T}}\left\vert f\right\vert ^{2}=\lambda \left\Vert f\right\Vert
^{2}\approx \left\Vert f\right\Vert ^{2}\mbox{, \ \ \ if \ }\lambda \approx 1%
\mbox{.}  \label{concentrated}
\end{equation}%
Thus, to count the number of degrees of freedom inside the region $T\times
\Omega $ for large $T$, Landau and Pollak \cite{LP} obtained the following
asymptotic estimate for the number of eigenvalues $\lambda _{n}$ of $%
P_{T,\Omega }$ which are close to $1$. For any $\delta >0$, 
\begin{equation}
\#\{n:\lambda _{n}>1-\delta \}\simeq \left\vert T\right\vert \left\vert
\Omega \right\vert /2\pi +C_{\delta }\log \left( \left\vert T\right\vert
\left\vert \Omega \right\vert \right) \text{,}  \label{estimate}
\end{equation}%
as $T\rightarrow \infty $, where $C_{\delta }$ is a constant depending only on $\delta$. The independence of $\delta $ allows to evaluate
the limit\ 
\begin{equation}
\lim_{r\rightarrow \infty }\frac{\eta (rT,\Omega )}{r}=\frac{\left\vert
T\right\vert \left\vert \Omega \right\vert }{2\pi }\mbox{,}  \label{L-P}
\end{equation}%
where $\eta (rT,\Omega )$ is the number of prolate spheroidal functions
essentially supported in the time- and bandlimited region $rT\times \Omega $%
. Within mathematical signal analysis (see, for instance the discussion in 
\cite[pag. 23]{Dau} and the recent book \cite{HogLak}), (\ref{estimate}) is
viewed as a mathematical formulation of the Nyquist rate, the fact that a
time- and bandlimited region $T\times \Omega $ corresponds to $\left\vert
T\right\vert \left\vert \Omega \right\vert /2\pi $ degrees of freedom.

\subsection{Slepian's solution: approximated dimension theorem}

With a view to solving the bandwith paradox, Slepian replaced the notions of
bandlimited and timelimited by more quantitative concepts, regarding signals
as $\epsilon $-timelimited in $T$ if the energy of the signal outside $T$ is
less than $\epsilon $ and $\epsilon $-bandlimited in $\Omega $ if the if the
energy of the Fourier transform of the signal outside $\Omega $ is less than 
$\epsilon $. Slepian associates $\epsilon $ to the precision of measuring
instruments and defines a flexible notion of $\epsilon $-approximate
dimension as follows. The set $F$ of signals is said to have approximate
dimension $N$ at level $\epsilon $ in the set $T$ if, for every $r\in F$,
there exist $a_{1},...a_{N}$ such that%
\begin{equation}
\int_{T}\left[ r(t)-\sum_{1}^{N}a_{j}g_{j}(t)\right] ^{2}dt<\epsilon 
\label{dimension}
\end{equation}%
and there is no set of $N-1$ functions that approximates every element of $F$
in such a way. Slepian%
\'{}%
s dimension theorem states that the approximated dimension $N(\Omega
,rT,\epsilon ,\epsilon 
{\acute{}}%
)$ at level $\epsilon 
{\acute{}}%
>\epsilon $ of the set $F_{\epsilon }$ of $\epsilon $-band and timelimited
functions, in the sense that $\left\Vert D_{rT}f-f\right\Vert \leq \epsilon $
and $\left\Vert B_{\Omega }f-f\right\Vert \leq \epsilon $, satisfies%
\begin{equation}
\lim_{r\rightarrow \infty }\frac{N(rT,\Omega ,\epsilon ,\epsilon 
{\acute{}}%
)}{r}=\frac{\left\vert T\right\vert \left\vert \Omega \right\vert }{2\pi }%
\mbox{.}  \label{Slepian}
\end{equation}%
Slepian's proof is also constructive. He defines a sequence of functions
using the prolates and their associated eigenvalues as follows:

$$
g_{j}(t)=\sqrt{\frac{\epsilon }{1-\lambda _{j}}}\phi _{j}(t)+\sqrt{\frac{%
\epsilon }{\lambda _{j}(1-\lambda _{j})}}\mathbf{1}_{\left[ -1,1\right]
}\left( \frac{2t}{T}\right) \phi _{j}(t)\mbox{.}
$$%
The $g_{j}$ are not complete in $F_{\epsilon }$, but Slepian has proved that
they are the best sequence to use for approximating functions in $%
F_{\epsilon }$.

\subsection{Pseudospectra enters the picture}

We remark that in the dimension theorems of Landau-Pollak and of Slepian,
the amount of energy outside $T\times \Omega $ does not appear in the
asymptotic limits (\ref{L-P}) and (\ref{Slepian}). With the aim of
developing an approximation theory of almost band-limited functions where
the number of degrees of freedom adjusts to the energy left outside $T\times
\Omega $, we have introduced a new sequence of functions which we call 
\emph{\ pseudo prolate spheroidal functions}. Our research program is not yet 
fully completed, but it is reasonable to expect these functions to have good
linear approximation properties of essentally band-limited functions, like
those recently proved for other orthogonal systems in (\cite{Jam}, \cite{Levitina}, \cite{Moumni}). Moreover,
we also expect the increase of the degrees of freedom to have
sparsity-promoting properties similar to the frame-based representations,
following the intuition provided by the "dictionary
example"\ \cite{Dictionary}: " The larger
and richer is my dictionary the shorter are the phrases I
compose."

We start by reformulating Landau-Pollack's approach in the following way:
instead of counting the number of eigenfunctions $f$ satisfying $P_{T,\Omega
}f=\lambda f$ which are associated with $\lambda \approx 1$, we count the
number of orthogonal functions such that $P_{T,\Omega }f\approx f$, in the
sense that the $L_{2}$ distance between $P_{T,\Omega }f\ $and$\ f$ is
smaller than a prescribed amount of energy $\epsilon $. Precisely, we assume 
$\left\Vert f\right\Vert =1$ and require that 
\begin{equation}
\left\Vert P_{T,\Omega }f-f\right\Vert ^{2}\leq \epsilon \mbox{.}  \label{b}
\end{equation}%
This measures the concentration of $f$ because $\epsilon $\ controls the
maximum amount of energy left outside $T\times \Omega $. For instance, if $%
\left\Vert D_{T}f-f\right\Vert \leq \epsilon $ and $\left\Vert B_{\Omega
}f-f\right\Vert \leq \epsilon $, then $\left\Vert P_{T,\Omega
}f-f\right\Vert ^{2}\leq 4\epsilon ^{2}$. The idea has been introduced in 
\cite{Measures}. It is based on the concept of \emph{pseudospectra of linear
operators}, which has found remarkable applications in the last decade (see 
\cite{Hansen}, \cite{Hansen2}, \cite{CPAM2004}, the surveys \cite{Tref97}
and \cite{SJ} or the book \cite{TE}). In general, $\lambda $ is an $\epsilon 
$-pseudoeigenvalue of $L$ if there exists $f$ with $\left\Vert f\right\Vert
=1$ such that $\left\Vert Lf-\lambda f\right\Vert \leq \epsilon $. Then $f
$ is said to be an $\epsilon $-pseudoeigenfunction corresponding to $\lambda $. 

As in the previous approaches, the set of $\epsilon $-localized functions in 
$T\times \Omega $ is not a linear space, making no sense to strictly talk
about its dimension. However, we can count the maximal number $\eta
(\epsilon ,rT,\Omega )$ of orthogonal functions satisfying (\ref{b}). In 
\cite{Measures}, using an explicit construction, we have shown that, as $%
r\rightarrow \infty $, the following inequalities hold:%
$$
\frac{\left\vert T\right\vert \left\vert \Omega \right\vert }{2\pi }%
(1+\epsilon )\leq \lim_{r\rightarrow \infty }\frac{\eta (\epsilon ,rT,\Omega
)}{r^{d}}\leq \frac{\left\vert T\right\vert \left\vert \Omega \right\vert }{%
2\pi }\left( 1-2\epsilon \right) ^{-1}\mbox{.}
$$%
Recently, we have obtained the sharp version of these inequalities: 
\begin{equation}
\lim_{r\rightarrow \infty }\frac{\eta (\epsilon ,rT,\Omega )}{r}=\left(
1-\epsilon \right) ^{-1}\frac{\left\vert T\right\vert \left\vert \Omega
\right\vert }{2\pi }\mbox{.}  \label{sharp}
\end{equation}%
Our proof of the lower inequality in (\ref{sharp}) is also constructive. As
in Slepian's approach and as in \cite{Measures}, the orthogonal functions
will be built in terms of the prolate spheroidal functions. We will describe
the construction of the functions achieving the sharp result (\ref{sharp}).
Since they result from a pseudospectra analogue of the spectral problem
defining the prolate spheroidal functions, we will call the corresponding
pseudoeigenfunctions\emph{\ pseudo prolate spheroidal functions}. \ Full
proofs of (\ref{sharp}) and other results will appear in \cite{Measures2}.

\section{Construction of the pseudo prolate spheroidal functions.}

\subsection{Time- and band- limiting operators}

A description of the general set-up of \cite{Landau1967} and \cite%
{Landau1975} follows. The sets $T$ and $\Omega $ are general subsets of
finite measure of $%
\mathbb{R}
^{d}$. Let%
$$
Ff(\xi )=\frac{1}{(2\pi )^{d/2}}\int_{\mathbb{R}^{d}}f(t)e^{-i\xi t}dt
$$%
denote the Fourier transform of a function $f\in L^{1}({\mathbb{R}}^{d})\cap
L^{2}({\mathbb{R}}^{d})$. The subspaces of $L^{2}({\mathbb{R}}^{d})$
consisting, respectively, of the functions supported in $T$ and of those
whose Fourier transform is supported in $\Omega $ are%
\begin{eqnarray*}
\mathcal{D}(T) &=&\{f\in L^{2}({\mathbb{R}}^{d}):f(x)=0,x\notin T\} \\
\mathcal{B}(\Omega ) &=&\{f\in L^{2}({\mathbb{R}}^{d}):Ff(\xi )=0,\xi \notin
\Omega \}.
\end{eqnarray*}%
Let $D_{T}$ be the orthogonal projection of $L^{2}({\mathbb{R}}^{d})$ onto $%
\mathcal{D}(T)$, given explicitly by the multiplication of a characteristic
function of the set $T$ by $f$: 
$$
D_{T}f(t)=\chi _{T}(t)f(t)
$$%
and let $B_{\Omega }$ be the orthogonal projection of $L^{2}({\mathbb{R}}%
^{d})$ onto $\mathcal{B}(\Omega )$, given explicitly as%
$$
B_{\Omega }f=F^{-1}\chi _{\Omega }Ff=\frac{1}{(2\pi )^{d/2}}\int_{{\mathbb{R}%
}^{d}}h(x-y)f(y)dy,
$$%
where $Fh=\chi _{\Omega }$.\ The following Theorem, comprising Lemma 1 and
Theorem 1 of \cite{Landau1975} gives important information concerning the
spectral problem associated with the operator $D_{rT}B_{\Omega }D_{rT}$. The
notation $o(r^{d})$ refers to behavior as $r\rightarrow \infty $.

\mbox{\textbf{Theorem A}}\cite{Landau1975}. \emph{The operator }$D_{rT}B_{\Omega
}D_{rT}$\emph{\ is bounded by }$1$\emph{, self-adjoint, positive, and
completely continuous. Denoting its set of eigenvalues, arranged in
nonincreasing order, by }$\{\lambda _{k}(r,T,\Omega )\}$\emph{, we have}

\begin{eqnarray*}
\sum_{k=0}^{\infty }\lambda _{k}(r,T,\Omega ) &=&r^{d}\left( 2\pi \right)
^{-d}\left\vert T\right\vert \left\vert \Omega \right\vert \\
\sum_{k=0}^{\infty }\lambda _{k}^{2}(r,T,\Omega ) &=&r^{d}\left( 2\pi
\right) ^{-d}\left\vert T\right\vert \left\vert \Omega \right\vert -o(r^{d})%
\mbox{.}
\end{eqnarray*}%
\emph{Moreover, given }$0<\gamma <1$\emph{, the number }$M_{r}(\gamma )$%
\emph{\ of eigenvalues which are not smaller than }$\gamma $\emph{,
satisfies, as }$r\rightarrow \infty $\emph{, }%
$$
M_{r}(\gamma )=\left( 2\pi \right) ^{-d}\left\vert T\right\vert \left\vert
\Omega \right\vert r^{d}+o(r^{d})\mbox{.}
$$

\subsection{Construction of the pseudo prolate spheroidal functions}

Suppose (\ref{b}) holds for a positive real $\varepsilon $. Let $\sigma >0$
be such that $\sigma ^{2}\leq \varepsilon $ and let $\mathcal{F}=\{\phi
_{k}\}$ be the normalized system of eigenfunctions (in the one dimension
interval case they are the prolates) of the operator $P_{rT,\Omega }$ with
eigenvalues $\lambda _{k}>1-\sigma $. Now, given $f\in L^{2}({\mathbb{R}}%
^{d})$, write%
\begin{equation}
f=\sum a_{k}\phi _{k}+h\mbox{,}  \label{expansf}
\end{equation}%
with $h\in $Ker$\left( P_{rT,\Omega }\right) $. Then%
\begin{equation}
P_{rT,\Omega }f=\sum a_{k}\lambda _{k}\phi _{k}  \label{fkern}
\end{equation}%
and%
\begin{eqnarray}
\left\Vert P_{rT,\Omega }f-f\right\Vert ^{2} &=&\left\Vert \sum (1-\lambda
_{k})a_{k}\phi _{k}+h\right\Vert ^{2}  \nonumber \\
&\leq &\sigma ^{2}\sum \left\vert a_{k}\right\vert ^{2}+\left\Vert
h\right\Vert ^{2}  \nonumber \\
&=&\sigma ^{2}\left\Vert f\right\Vert ^{2}+(1-\sigma ^{2})\left\Vert
h\right\Vert ^{2}\mbox{.}  \label{third}
\end{eqnarray}%
For the given $\sigma >0$ we pick a real number $\gamma $ such that 
\begin{equation}
\sigma ^{2}+(1-\sigma ^{2})\gamma =\varepsilon \mbox{.}  \label{epsilon}
\end{equation}%
Writing this as $\gamma =(\varepsilon -\sigma ^{2})/(1-\sigma ^{2})$ it's
clear that $\gamma $ is a positive increasing function of $\sigma $, and $%
\gamma \rightarrow \varepsilon $ as $\sigma \rightarrow 0$. Now take $n=\#%
\mathcal{F}$, define $\Gamma =\left\langle \mathcal{F}\right\rangle $ and
let $m$ be a positive integer (its value will be made precise later). Choose 
$h_{1},h_{2},\dots ,h_{m}$ orthonormal functions in Ker$\left( P_{rT,\Omega
}\right) $, and let $\Lambda $ be the space spanned by this functions. This
can be done since Ker$\left( P_{rT,\Omega }\right) $ has infinite dimension,
due to the inclusion $\mathcal{D}({\mathbb{R}}^{d}-rT)\subset $Ker$\left(
P_{rT,\Omega }\right) $. Note that this $m$ functions together with the $n$
functions of $\mathcal{F}$ form a orthonormal basis of $\Gamma \oplus
\Lambda $, since the first are orthogonal to the latter. We now \emph{define
the pseudoeigenfunctions} as a second orthonormal basis of $\Gamma \oplus
\Lambda $, denoted by $\{\Phi _{j}\}_{j=1}^{m+n}$, with 
\begin{equation}
\Phi _{j}=\psi _{j}+\rho _{j},  \label{pseudocandidates}
\end{equation}%
$\psi _{j}\in \Gamma $ and $\rho _{j}\in \Lambda $ for $j\in \{1,2,\dots
,m+n\}$. The proof of the lower inequality in (\ref{sharp}) requires the
construction of the functions (\ref{pseudocandidates}) in such a way that 
\begin{equation}
\Vert \rho _{j}\Vert ^{2}=\frac{m}{m+n},\quad j=1,2,\dots ,m+n\mbox{.}
\label{disterrcond}
\end{equation}%
This will be done using a linear algebra argument detailed in the next
paragraph.

Consider the automorphism $Q$ in $\Gamma \oplus \Lambda $ that maps the
first basis to the functions in (\ref{pseudocandidates}). One can see $Q$ as
an orthogonal $(m+n)\times (m+n)$ matrix of the form%
$$
Q=[Q^{\Gamma }\;Q^{\Lambda }]_{(m+n)\times (m+n)}\mbox{,}
$$%
where its first $n$ columns $Q^{\Gamma }$ map $\Gamma $ to $\Gamma \oplus
\Lambda $ and the last $m$ columns $Q^{\Lambda }$ map $\Lambda $ to $\Gamma
\oplus \Lambda $. Then, since $\Vert h_{j}\Vert =1$, the condition (\ref%
{disterrcond}) is equivalent to 
\begin{equation}
\Vert Q_{j}^{\Lambda }\Vert ^{2}=\frac{m}{m+n},\quad j=1,2,\dots ,m+n,
\label{disterrcondmat}
\end{equation}%
where $Q_{j}^{\Lambda }$ denotes the $j$th line of $Q^{\Lambda }$. Let $X$
be the Discrete Fourier Transform matrix of order $m+n$, with entries 
$$
X_{ij}=\frac{1}{\sqrt{m+n}}\omega ^{ij},\quad i,j=0,1,\dots ,m+n-1,
$$%
where $\omega =e^{-\frac{2\pi i}{m+n}}$ is the $(m+n)$th-root of the unity.
Then define the $(m+n)\times (m+n)$ matrix $X^{\prime }$ as 
\begin{equation}
X_{ij}^{\prime }=\Re (X_{ij})+\Im (X_{ij}),\quad i,j=0,1,\dots ,m+n-1,
\end{equation}%
One can check that this matrix is orthogonal (more details will be given in \cite{Measures2}).
Now we finally define $Q$ as a permutation of the columns of $X^{\prime }$,
depending on the parity of $m$. If $m$ is even, we choose the last $m$
columns of $Q$ to be the $1$ to $m/2$ and the last $m/2$ columns of $%
X^{\prime }$. This leads to 
\begin{eqnarray*}
\Vert Q_{j}^{\Lambda }\Vert ^{2} &=&\sum_{k=n+1}^{m+n}Q_{jk}^{2} \\
&=&\frac{1}{m+n}\sum_{k=1}^{m/2}(a_{jk}+b_{jk})^{2}+(a_{-jk}+b_{-jk})^{2} \\
&=&\frac{1}{m+n}\sum_{k=1}^{m/2}(a_{jk}+b_{jk})^{2}+(a_{jk}-b_{jk})^{2} \\
&=&\frac{1}{m+n}\sum_{k=1}^{m/2}2a_{jk}^{2}+2b_{jk}^{2} \\
&=&\frac{m}{m+n},
\end{eqnarray*}%
since $a_{jk}^{2}+b_{jk}^{2}=1$, thus $Q$ satisfies (\ref{disterrcondmat}).
If $m$ is odd, we add to this columns the column $0$, which has all entries
equal to $1/\sqrt{m+n}$, the additional calculations in this case are
trivial. We have finally proved that there are $m+n$ functions as in (\ref%
{pseudocandidates}) which verify (\ref{disterrcond}). Since $\psi _{j}$ are
linear combinations of elements of $\mathcal{F}=\{\phi _{k}\}$, and $\rho
_{j}\in \Ker\left( P_{rT,\Omega }\right) $, (\ref{pseudocandidates}) is a
representation of the form (\ref{expansf}). We can now apply (\ref{third})
and (\ref{disterrcond}) to obtain%
\begin{eqnarray}
\left\Vert P_{rT,\Omega }\Phi _{j}-\Phi _{j}\right\Vert ^{2} &\leq &\sigma
^{2}\left\Vert \Phi _{j}\right\Vert ^{2}+(1-\sigma ^{2})\Vert \rho _{j}\Vert
^{2} \nonumber \\
&=&\sigma ^{2}+(1-\sigma ^{2})\frac{m}{m+n}  \label{mchoose}
\end{eqnarray}%
We now choose $m$ so that (\ref{mchoose}) is at most $\varepsilon $, or
equivalently, $\frac{m}{m+n}\leq \gamma $. Clearly, this happens if and only
if $m\leq n\gamma /(1-\gamma )$. Choosing the biggest $m$ which verifies
this condition, leads to $m\geq \frac{n\gamma }{1-\gamma }-1$. We now use
Theorem A (the fact that $n=\#\mathcal{F}=r^{d}\left( 2\pi \right)
^{-d}\left\vert T\right\vert \left\vert \Omega \right\vert +o(r^{d})$) and
this last inequality 
\begin{eqnarray*}
\#\{\Phi _{j}\}_{j=1}^{m+n} &=&m+n \\
&\geq &n\left( \frac{\gamma }{1-\gamma }+1\right) -1 \\
&=&\left( \frac{1}{1-\gamma }\right) n-1 \\
&=&({1-\gamma })^{-1}\left( r^{d}(2\pi )^{-d}|T||\Omega |+o(r^{d})\right) -1
\\
&=&({1-\gamma })^{-1}\left( r^{d}(2\pi )^{-d}|T||\Omega |+o(r^{d})\right) ,
\end{eqnarray*}%
since $1=o(r^{d})$. We have obtained by construction the pseudo prolate
spheroidal functions $\{\Phi _{j}\}_{j=1}^{m+n}$. They are also orthonormal
and verify (\ref{b}).

The lower inequality in (\ref{sharp}) is now a simple consequence of this
construction. Denote by $M^{-}(rT,\Omega ,\epsilon )$ the minimum number of
orthonormal functions satisfying (\ref{b}). Then, 
\begin{eqnarray*}
M^{-}(rT,\Omega ,\epsilon ) &\geq& \#\left[ \cup _{i=1}^{l}\{\Phi
_{j}^{(i)}\}_{j=1}^{n+1}\right] \\
& \geq& ({1-\gamma })^{-1}r^{d}\left( 2\pi
\right) ^{-d}\left\vert T\right\vert \left\vert \Omega \right\vert +o(r^{d})%
\mbox{.}
\end{eqnarray*}
Finally, take $\sigma \rightarrow 0$, so that $\gamma \rightarrow \epsilon $
to yield%
\begin{eqnarray*}
M^{-}(rT,\Omega ,\epsilon )&\geq& \#\left[ \cup _{i=1}^{l}\{\Phi
_{j}^{(i)}\}_{j=1}^{n+1}\right] \\
&\geq& (1-\epsilon )^{-1}r^{d}\left( 2\pi
\right) ^{-d}\left\vert T\right\vert \left\vert \Omega \right\vert +o(r^{d})%
\mbox{.}
\end{eqnarray*}

\begin{remark}
The above construction applies to several settings where properties similar
to Theorem A are available, as in \cite{AB} and \cite{FeiNow},\cite{AGR}.
\end{remark}

\end{document}